\documentclass[12pt]{article}%
\usepackage{amsmath}
\usepackage{amsfonts}
\usepackage{amssymb}
\usepackage{graphicx}%
\newtheorem{theorem}{Theorem}

\newtheorem{corollary}[theorem]{Corollary}

\newtheorem{example}[theorem]{Example}

\newtheorem{lemma}[theorem]{Lemma}

\newtheorem{remark}[theorem]{Remark}

\newenvironment{proof}[1][Proof]{\noindent\textbf{#1.} }{\ \rule{0.5em}{0.5em}}
\begin{document}

\date{}
\title{Asymptotic analysis of powers of matrices}
\author{Diego Dominici \thanks{e-mail: dominicd@newpaltz.edu}\\Department of Mathematics\\State University of New York at New Paltz\\75 S. Manheim Blvd. Suite 9\\New Paltz, NY 12561-2443\\USA\\Phone: (845) 257-2607\\Fax: (845) 257-3571 }
\maketitle

\begin{abstract}
We analyze the representation of $A^{n}$ as a linear combination of
$A^{j},\ 0\leq j\leq k-1,$ where $A$ is a $k\times k$ matrix. We obtain a
first order asymptotic approximation of $A^{n}$ as $n\rightarrow\infty,$
without imposing any special conditions on $A.$

We give some examples showing the application of our results.

\end{abstract}

Keywords: Matrix powers, asymptotic approximations of integrals, generating functions. 

MSC-class: 15A15 (Primary) 34E05 (Secondary)

\section{Introduction}

In a recent article \cite{MR2125272}, Abu-Saris and Ahmad showed how to
compute the powers of a matrix without having to compute its eigenvalues.
Their main result was:

\begin{theorem}
If $A$ is a $k\times k$ matrix with characteristic polynomial%
\begin{equation}
P(x)=x^{k}+%
{\displaystyle\sum\limits_{j=0}^{k-1}}
a_{j}x^{j}, \label{p}%
\end{equation}
then,%
\begin{equation}
A^{n}=%
{\displaystyle\sum\limits_{j=0}^{k-1}}
b_{j}(n)A^{j},\quad n\geq k \label{powers}%
\end{equation}
where%
\begin{align}
b_{j}(k)  &  =-a_{j},\quad0\leq j\leq k-1,\quad b_{-1}(n)=0,\quad n\geq
k,\label{req}\\
b_{j}(n+1)  &  =b_{j-1}(n)-a_{j}b_{k-1}(n),\quad n\geq k,\quad0\leq j\leq
k-1.\nonumber
\end{align}

\end{theorem}

The purpose of this paper is to find an asymptotic representation for the
numbers $b_{j}(n)$ as $n\rightarrow\infty,$ which using (\ref{powers}) will
give an asymptotic representation of $A^{n}$ for large $n.$ Since the
coefficients $b_{j}(n)$ depend only on $P(x),$ our estimates will be valid for
similar matrices.

The asymptotic behavior of powers of matrices has been considered before by
other authors. In \cite{MR0056568} and \cite{MR0063337}, Gautschi computed
upper bounds for $A^{n}$ and $\left\Vert A^{n}\right\Vert ,$ where $\left\Vert
A\right\Vert $ is a norm of $A.$ Estimates of $\left\Vert A^{n}\right\Vert $
were also studied in \cite{MR1867596}, \cite{MR1803180}, \cite{MR806832},
\cite{MR484965} and \cite{MR520625}.

In \cite{MR578322}, Friedland and Schneider considered the matrix%
\[
B^{(m)}=A^{m}\left(  I+\cdots+A^{q-1}\right)  ,\quad m\geq1
\]
where $A$ is a nonnegative matrix and $q$ is a certain positive integer. They
proved a theorem on the growth of $B^{(m)}$ under the assumption that the
spectral radius of $A$ is equal to one. Powers of nonnegative matrices were
also analyzed by Lindqvist in \cite{MR986894}. \ Rothblum \cite{MR618637},
obtained Cesaro asymptotic expansions of $%
{\displaystyle\sum\limits_{i=0}^{N}}
A^{i},$ where $A$ is a complex matrix with spectral radius less than or equal
to one.

This paper is organized as follows: In Section \ref{section1} we find an
integral representation for the exponential generating function $G_{j}(z)$ of
the coefficients $b_{j}(n).$ We obtain exact formulas for $G_{j}(z)$ and
$b_{j}(n)$ in the special case of the matrix $A$ having $k$ distinct
eigenvalues. We conclude the section with some examples.

In Section \ref{section2} we give an exact representation and a first order
asymptotic approximation for $b_{j}(n),$ as $n\rightarrow\infty$. We consider
the cases of simple and multiple eigenvalues. Our formulas are relatively easy
to implement and offer very accurate estimates of $b_{j}(n),$ and therefore of
$A^{n},$ for large $n.$ We present some examples for different cases of
$P(x).$

\section{ Generating function}

\label{section1}

In this section we shall find and exponential generating function for the
coefficients $b_{j}(n).$ First, let us define the spectral radius $\rho(A)$ of
the matrix $A$ by%
\begin{equation}
\rho(A)=\max\left\{  \left\vert \lambda\right\vert \ |\ P(\lambda)=0\right\}
.\label{rho}%
\end{equation}

\begin{theorem}
Let $G_{j}(z)$ be defined by
\begin{equation}
G_{j}(z)=%
{\displaystyle\sum\limits_{n\geq0}}
b_{j}(n+k)\frac{z^{n}}{n!}, \label{G}%
\end{equation}
Then, we have%
\begin{equation}
G_{j}(z)=-\frac{1}{2\pi\mathrm{i}}%
{\displaystyle\int\limits_{c-\mathrm{i}\infty}^{c+\mathrm{i}\infty}}
s^{k-j-1}\frac{p_{j}(s)}{P(s)}e^{zs}ds, \label{G1}%
\end{equation}
where $c>\rho(A)$,%
\begin{equation}
p_{j}(s)=%
{\displaystyle\sum\limits_{l=0}^{j}}
a_{l}x^{l},\quad0\leq j\leq k-1, \label{Pr}%
\end{equation}
and $P(s)$ is the characteristic polynomial of $A$ defined in (\ref{Pr}).
\end{theorem}

\begin{proof}
If we use (\ref{G}) in (\ref{req}), we obtain%
\begin{equation}
G_{j}^{\prime}=G_{j-1}-a_{j}G_{k-1},\quad G_{j}(0)=-a_{j},\quad0\leq j\leq k-1
\label{diffG}%
\end{equation}
with $G_{-1}(z)=0.$ Taking the Laplace transform of $G_{j}(z),$%
\[
L_{j}(s)=%
{\displaystyle\int\limits_{0}^{\infty}}
G_{j}(z)e^{-zs}dz
\]
in (\ref{diffG}) we get
\begin{equation}
sL_{j}+a_{j}=L_{j-1}-a_{j}L_{k-1},\quad0\leq j\leq k-1 \label{sys}%
\end{equation}
and $L_{-1}(z)=0.$
\end{proof}

The solution of (\ref{sys}) is given by%
\begin{equation}
L_{j}(s)=-s^{k-j-1}\frac{p_{j}(s)}{P(s)},\quad0\leq j\leq k-1 \label{L}%
\end{equation}
where%
\[
p_{j}(s)=%
{\displaystyle\sum\limits_{l=0}^{j}}
a_{l}x^{l},\quad0\leq j\leq k-1,
\]
and $p_{-1}(s)=0$. Inverting the Laplace transform in (\ref{L}) the theorem follows.

\begin{remark}
Since
\[
\underset{\left\vert s\right\vert \rightarrow\infty}{\lim}L_{j}(s)=0,\quad
0\leq j\leq k-1
\]
we can replace the Bromwich contour in (\ref{G1}) with a circle $\mathcal{C}$
of radius $R$ centered at the origin \cite{MR616824}%
\begin{equation}
G_{j}(z)=-\frac{1}{2\pi\mathrm{i}}%
{\displaystyle\int\limits_{\mathcal{C}}}
s^{k-j-1}\frac{p_{j}(s)}{P(s)}e^{zs}ds \label{GC}%
\end{equation}
with $R>\rho(A)$.
\end{remark}

\begin{corollary}
If
\[
P(s)=(s-\lambda_{1})(s-\lambda_{2})\cdots(s-\lambda_{k}),
\]
where the eigenvalues $\lambda_{i}$ are all distinct, then%
\begin{equation}
G_{j}(z)=-%
{\displaystyle\sum\limits_{l=1}^{k}}
\left(  \lambda_{l}\right)  ^{k-j-1}\frac{p_{j}(\lambda_{l})}{P^{\prime
}(\lambda_{l})}\exp\left(  \lambda_{l}z\right)  \label{G2}%
\end{equation}
and%
\begin{equation}
b_{j}(n)=-%
{\displaystyle\sum\limits_{l=1}^{k}}
\left(  \lambda_{l}\right)  ^{n-j-1}\frac{p_{j}(\lambda_{l})}{P^{\prime
}(\lambda_{l})}. \label{b1}%
\end{equation}

\end{corollary}

\begin{proof}
Applying the residue theorem to (\ref{GC}) we obtain%
\begin{equation}
G_{j}(z)=-%
{\displaystyle\sum\limits_{P\left(  \lambda\right)  =0}}
\operatorname{Re}\mathrm{s}\left[  s^{k-j-1}\frac{p_{j}(s)}{P(s)}%
e^{zs};\lambda\right]  \label{Gres}%
\end{equation}
which in turn gives (\ref{G2}) after computing%
\begin{gather*}
\operatorname{Re}\mathrm{s}\left[  s^{k-j-1}\frac{p_{j}(s)}{P(s)}%
e^{zs};\lambda_{l}\right]  =\underset{s\rightarrow\lambda_{l}}{\lim}%
s^{k-j-1}p_{j}(s)e^{zs}\frac{(s-\lambda_{l})}{P(s)}\\
=\underset{s\rightarrow\lambda_{l}}{\lim}s^{k-j-1}p_{j}(s)e^{zs}\frac
{1}{P^{\prime}(s)}=\left(  \lambda_{l}\right)  ^{k-j-1}\frac{p_{j}(\lambda
_{l})}{P^{\prime}(\lambda_{l})}\exp\left(  \lambda_{l}z\right)  .
\end{gather*}

Writing (\ref{G2}) as%
\[
G_{j}(z)=-%
{\displaystyle\sum\limits_{l=1}^{k}}
\left(  \lambda_{l}\right)  ^{k-j-1}\frac{p_{j}(\lambda_{l})}{P^{\prime
}(\lambda_{l})}%
{\displaystyle\sum\limits_{n\geq0}}
\left(  \lambda_{l}\right)  ^{n}\frac{z^{n}}{n!}%
\]
and changing the order of summation, we have%
\[
G_{j}(z)=%
{\displaystyle\sum\limits_{n\geq0}}
\left[  -%
{\displaystyle\sum\limits_{l=1}^{k}}
\left(  \lambda_{l}\right)  ^{n+k-j-1}\frac{P_{j}(\lambda_{l})}{P^{\prime
}(\lambda_{l})}\right]  \frac{z^{n}}{n!}%
\]
which implies%
\[
b_{j}\left(  n+k\right)  =-%
{\displaystyle\sum\limits_{l=1}^{k}}
\left(  \lambda_{l}\right)  ^{n+k-j-1}\frac{p_{j}(\lambda_{l})}{P^{\prime
}(\lambda_{l})}.
\]

\end{proof}

\begin{example}
\label{ex1} In \cite{MR2125272} the authors considered the following examples:

\begin{enumerate}
\item
\[
P(x)=x^{3}-7x^{2}+16x-12=\left(  x-2\right)  ^{2}\left(  x-3\right)  .
\]
Using (\ref{L}) we have%
\begin{align*}
L_{0}(s)  &  =\frac{-s^{2}\left(  12\right)  }{\left(  s-2\right)  ^{2}\left(
s-3\right)  }\\
L_{1}(s)  &  =\frac{-s\left(  16s-12\right)  }{\left(  s-2\right)  ^{2}\left(
s-3\right)  }\\
L_{2}(s)  &  =\frac{-\left(  -7x+16x-12\right)  }{\left(  s-2\right)
^{2}\left(  s-3\right)  }%
\end{align*}
and inverting we obtain%
\begin{align*}
G_{0}(z)  &  =-12\left(  8+4z\right)  e^{2z}+108e^{3z}\\
G_{1}(z)  &  =4\left(  23+10z\right)  e^{2z}-108e^{3z}\\
G_{2}(z)  &  =-4\left(  5+2z\right)  e^{2z}+27e^{3z}.
\end{align*}
Expanding in series we get%
\[
G_{0}(z)=-96%
{\displaystyle\sum\limits_{n\geq0}}
2^{n}\frac{z^{n}}{n!}-48%
{\displaystyle\sum\limits_{n\geq0}}
2^{n-1}n\frac{z^{n}}{n!}+108%
{\displaystyle\sum\limits_{n\geq0}}
3^{n}\frac{z^{n}}{n!}%
\]
and from (\ref{G}) we conclude that%
\[
b_{0}\left(  n+3\right)  =-96\times2^{n}-48\times2^{n-1}n+108\times3^{n}%
\]
or%
\begin{align*}
b_{0}\left(  n\right)   &  =-96\times2^{n-3}-48\times2^{n-4}\left(
n-3\right)  +108\times3^{n-3}\\
&  =-3\left(  1+n\right)  \times2^{n}+4\times3^{n}.
\end{align*}
Similar calculations give%
\begin{align*}
b_{1}\left(  n\right)   &  =\left(  4+\frac{5}{2}n\right)  \times2^{n}%
-4\times3^{n}\\
b_{2}\left(  n\right)   &  =-\left(  1+\frac{1}{2}n\right)  \times2^{n}+3^{n},
\end{align*}
in agreement with the results shown in \cite{MR2125272}.

\item
\[
P(x)=x^{3}-5x^{2}+6x=x\left(  x-2\right)  \left(  x-3\right)  .
\]
We can apply (\ref{b1}) directly and obtain%
\begin{align*}
b_{0}\left(  n\right)   &  =0\\
b_{1}\left(  n\right)   &  =\frac{3}{2}\times2^{n}-\frac{2}{3}\times3^{n}\\
b_{2}\left(  n\right)   &  =-\frac{1}{2}\times2^{n}+\frac{1}{3}\times3^{n}.
\end{align*}

\item
\[
P(x)=x^{5}-5x^{4}+10x^{3}-20x^{2}-15x-4=\left(  x-4\right)  \left(
x^{4}-x^{3}+6x^{2}+4x+1\right)  .
\]
Although (as the authors noted) MAPLE is unable to compute the zeros of $P(x)$
exactly, it can provide us with very accurate numerical approximations%
\begin{align*}
\lambda_{1}  &  =4\\
\lambda_{2}  &  =0.8090169944+2.489898285\mathrm{i}\\
\lambda_{3}  &  =0.8090169944-2.489898285\mathrm{i}\\
\lambda_{4}  &  =-0.3090169944+0.2245139883\mathrm{i}\\
\lambda_{5}  &  =-0.3090169944-0.2245139883\mathrm{i}%
\end{align*}
which we can use in (\ref{b1}) to get%
\[
b_{j}(n)=\frac{C_{j}}{305}4^{n}-%
{\displaystyle\sum\limits_{l=2}^{5}}
\left(  \lambda_{l}\right)  ^{n-j-1}\frac{p_{j}(\lambda_{l})}{P^{\prime
}(\lambda_{l})}%
\]
with
\[
C_{0}=1,C_{1}=4,C_{2}=6,C_{3}=-1,C_{4}=1.
\]
Note that, for $0\leq j\leq4,$ we have
\[%
{\displaystyle\sum\limits_{l=2}^{5}}
\left(  \lambda_{l}\right)  ^{n-j-1}\frac{p_{j}(\lambda_{l})}{P^{\prime
}(\lambda_{l})}=O\left(  \left\vert \lambda_{2}\right\vert ^{n}\right)
=O\left(  2.618^{n}\right)
\]
as $n\rightarrow\infty.$
\end{enumerate}
\end{example}

\begin{example}
Let $A$ be the matrix
\[
A=%
\begin{pmatrix}
1 & 2\\
-1 & -1
\end{pmatrix}
\]
with characteristic polynomial
\[
P(x)=x^{2}+1=(x-\mathrm{i})(x+\mathrm{i}).
\]
Using (\ref{b1}) we have%
\begin{align*}
b_{0}(n)  &  =\cos\left(  \frac{\pi}{2}n\right) \\
b_{1}(n)  &  =\sin\left(  \frac{\pi}{2}n\right)
\end{align*}
and from (\ref{powers}) we get%
\[
A^{n}=%
\begin{pmatrix}
\cos\left(  \frac{\pi}{2}n\right)  +\sin\left(  \frac{\pi}{2}n\right)  &
2\sin\left(  \frac{\pi}{2}n\right) \\
-\sin\left(  \frac{\pi}{2}n\right)  & \cos\left(  \frac{\pi}{2}n\right)
-\sin\left(  \frac{\pi}{2}n\right)
\end{pmatrix}
.
\]
In particular, we have%
\[
A^{n}=\left\{
\begin{array}
[c]{c}%
I,\quad n\equiv0(4)\\
A,\quad n\equiv1(4)\\
-I,\quad n\equiv2(4)\\
-A,\quad n\equiv3(4)
\end{array}
\right.
\]
where $I$ denotes the identity matrix.
\end{example}

\begin{example}
This example appeared in \cite{MR1739433}. Let $A$ be the matrix
\[
A=%
\begin{pmatrix}
1 & 1\\
1 & 0
\end{pmatrix}
\]
with characteristic polynomial
\[
P(x)=x^{2}-x-1=(x-\alpha)(x-\beta),
\]
where%
\[
\alpha=\frac{1}{2}\left(  1+\sqrt{5}\right)  ,\quad\beta=\frac{1}{2}\left(
1-\sqrt{5}\right)  .
\]
Then, from (\ref{b1}), we have%
\begin{align*}
b_{0}(n)  &  =\frac{1}{\sqrt{5}}\left(  \alpha^{n-1}-\beta^{n-1}\right)
=f_{n-1}\\
b_{1}(n)  &  =\frac{1}{\sqrt{5}}\left(  \alpha^{n}-\beta^{n}\right)  =f_{n}%
\end{align*}
where $f_{n}$ is the nth Fibonacci number. Thus,%
\[
A^{n}=%
\begin{pmatrix}
f_{n}+f_{n-1} & f_{n}\\
f_{n} & f_{n-1}%
\end{pmatrix}
=%
\begin{pmatrix}
f_{n+1} & f_{n}\\
f_{n} & f_{n-1}%
\end{pmatrix}
.
\]

\end{example}

\section{Asymptotic analysis}

\label{section2}

We begin by finding an integral representation of the coefficients $b_{j}(n).$

\begin{lemma}
The numbers $b_{j}(n)$ can be represented as%
\begin{equation}
b_{j}(n)=-\frac{1}{2\pi\mathrm{i}}%
{\displaystyle\int\limits_{\mathcal{C}}}
s^{n-j-1}\frac{p_{j}(s)}{P(s)}ds \label{binteg}%
\end{equation}
where $\mathcal{C}$ \ is a circle of radius $R>\rho(A)$ centered at the origin
and the polynomials $p_{j}(s)$ were defined in (\ref{Pr}).
\end{lemma}

\begin{proof}
Since the power series
\[
e^{zs}=%
{\displaystyle\sum\limits_{n\geq0}}
s^{n}\frac{z^{n}}{n!}%
\]
converges uniformly on $\left\vert s\right\vert \leq R,$ we can interchange
integration and summation in (\ref{GC}) and obtain%
\[
G_{j}(z)=%
{\displaystyle\sum\limits_{n\geq0}}
\left[  -\frac{1}{2\pi\mathrm{i}}%
{\displaystyle\int\limits_{\mathcal{C}}}
s^{n+k-j-1}\frac{p_{j}(s)}{P(s)}ds\right]  \frac{z^{n}}{n!}.
\]
Then, (\ref{G}) implies%
\[
b_{j}(k+n)=-\frac{1}{2\pi\mathrm{i}}%
{\displaystyle\int\limits_{\mathcal{C}}}
s^{n+k-j-1}\frac{p_{j}(s)}{P(s)}ds,\quad0\leq j\leq k-1
\]
and the result follows.
\end{proof}

\begin{remark}
An alternative method for approximating the coefficients $b_{j}(n)$ is to
write (\ref{binteg}) as%
\begin{equation}
b_{j}(n)=-\frac{R^{n-j}}{2\pi}%
{\displaystyle\int\limits_{0}^{2\pi}}
\exp\left[  \mathrm{i}t(n-j)\right]  \frac{p_{j}(Re^{\mathrm{i}t}%
)}{P(Re^{\mathrm{i}t})}dt \label{binteg2}%
\end{equation}
with $R>\rho(A)$ and to compute the integral (\ref{binteg2}) numerically. This
approach offers the advantage of avoiding the computation of the eigenvalues
of $A.$
\end{remark}

We now have all the necessary elements to establish our main theorem.

\begin{theorem}
Let
\[
\rho(A)=\left\vert \lambda\right\vert >\left\vert \lambda_{2}\right\vert
>\cdots>\left\vert \lambda_{r}\right\vert
\]
be the eigenvalues of the matrix $A,$ i.e.,%
\begin{equation}
P(x)=\left(  x-\lambda\right)  ^{m}\left(  x-\lambda_{2}\right)  ^{m_{2}%
}\cdots\left(  x-\lambda_{r}\right)  ^{m_{r}} \label{P1}%
\end{equation}
with $r\leq k.$ Then,
\begin{equation}
b_{j}(n)\sim-\lambda^{n-m-j}\frac{p_{j}(\lambda)}{P^{\left(  m\right)
}(\lambda)}m!\binom{n-k}{m-1},\quad n\rightarrow\infty, \label{basymp}%
\end{equation}
where%
\[
P^{\left(  m\right)  }(\lambda)=\left.  \frac{d^{m}P}{ds^{m}}\right\vert
_{s=\lambda}.
\]

\end{theorem}

\begin{proof}
To find an asymptotic approximation of (\ref{binteg}), we shall use a modified
version of Darboux's Method \cite{MR1851050}. We write
\[
s^{n-j-1}\frac{p_{j}(s)}{P(s)}=s^{n-k}\times\frac{s^{k-j-1}p_{j}(s)}{P(s)},
\]
so that $\deg\left(  s^{k-j-1}p_{j}\right)  =k-1$ and $\deg\left(  P\right)
=k.$

From (\ref{P1}) we have%
\begin{equation}
\frac{s^{k-j-1}p_{j}(s)}{P(s)}\sim\frac{\lambda^{k-j-1}p_{j}(\lambda)}{\left(
s-\lambda\right)  ^{m}g(\lambda)},\quad s\rightarrow\lambda\label{s1}%
\end{equation}
where%
\[
g(x)=\left(  x-\lambda_{2}\right)  ^{m_{2}}\cdots\left(  x-\lambda_{r}\right)
^{m_{r}}.
\]
Using the Binomial Theorem, we obtain%
\begin{equation}
s^{n-k}=%
{\displaystyle\sum\limits_{l=0}^{n-k}}
\left(  s-\lambda\right)  ^{l}\binom{n-k}{l}\lambda^{n-k-l}. \label{s2}%
\end{equation}
Combining (\ref{s1}) and (\ref{s2}), we get%
\[
s^{n-j-1}\frac{p_{j}(s)}{P(s)}\sim\lambda^{n-m-j}\frac{p_{j}(\lambda
)}{g(\lambda)}\binom{n-k}{m-1}\frac{1}{\left(  s-\lambda\right)  },\quad
s\rightarrow\lambda
\]
and therefore
\begin{equation}
b_{j}(n)=-\frac{1}{2\pi\mathrm{i}}%
{\displaystyle\int\limits_{\mathcal{C}}}
s^{n-j-1}\frac{p_{j}(s)}{P(s)}ds\sim-\lambda^{n-m-j}\frac{p_{j}(\lambda
)}{g(\lambda)}\binom{n-k}{m-1},\quad n\rightarrow\infty. \label{b11}%
\end{equation}

To find the value of $g(\lambda),$ we use L'Hopital's Theorem
\begin{equation}
g(\lambda)=\underset{s\rightarrow\lambda}{\lim}\frac{P(s)}{\left(
s-\lambda\right)  ^{m}}=\underset{s\rightarrow\lambda}{\lim}\frac{P^{(m)}%
(s)}{m!}=\frac{P^{(m)}(\lambda)}{m!}. \label{g1}%
\end{equation}
Replacing (\ref{g1}) in (\ref{b11}), we obtain (\ref{basymp}).
\end{proof}

\begin{remark}
Note that when $m=1$ we recover the leading term in (\ref{b1}).
\end{remark}

If more than one eigenvalue has absolute value equal to the spectral radius of
$A,$ the asymptotic behavior of $b_{j}(n)$ can be obtained by adding the
contributions from each eigenvalue. We state this formally in the following corollary.

\begin{corollary}
If
\[
\rho(A)=\left\vert \lambda_{1}\right\vert =\left\vert \lambda_{2}\right\vert
=\cdots=\left\vert \lambda_{r}\right\vert ,
\]
with respective multiplicities $m_{1},m_{2},\ldots,m_{r},$ then%
\begin{equation}
b_{j}(n)\sim-%
{\displaystyle\sum\limits_{l=1}^{r}}
\left(  \lambda_{l}\right)  ^{n-m_{l}-j}\frac{p_{j}(\lambda_{l})}{P^{\left(
m_{l}\right)  }(\lambda_{l})}\left(  m_{l}\right)  !\binom{n-k}{m_{l}-1},\quad
n\rightarrow\infty. \label{basymp2}%
\end{equation}

\end{corollary}

\begin{remark}
Since%
\[
\binom{n-k}{m_{l}-1}\sim\frac{1}{\left(  m_{l}-1\right)  !}n^{m_{l}-1},\quad
n\rightarrow\infty
\]
we have%
\[
b_{j}(n)\sim-%
{\displaystyle\sum\limits_{l=1}^{r}}
\left(  \lambda_{l}\right)  ^{n-m_{l}-j}\frac{p_{j}(\lambda_{l})}{P^{\left(
m_{l}\right)  }(\lambda_{l})}m_{l}n^{m_{l}-1},\quad n\rightarrow\infty.
\]
Therefore, in the case of several eigenvalues located on the circle
$\left\vert s\right\vert =\rho(A),$ the dominant term in (\ref{basymp2}) will
correspond to the eigenvalue with the greatest multiplicity.
\end{remark}

\begin{example}
In Example \ref{ex1} (1) we consider
\[
P(x)=x^{3}-7x^{2}+16x-12=\left(  x-2\right)  ^{2}\left(  x-3\right)  .
\]
In this case, $\lambda=2,$ $m=2$ and $k=3.$ From (\ref{basymp}), we get%
\begin{align*}
b_{0}(n) &  \sim4\times3^{n}\\
b_{1}(n) &  \sim-4\times3^{n}\\
b_{2}(n) &  \sim3^{n}%
\end{align*}
which are the leading terms in the solution previously obtained.
\end{example}

\begin{example}
We now consider the case of more than one eigenvalue having absolute value
equal to $\rho(A).$ Let%
\[
P(x)=x^{4}+x^{3}-15x^{2}-9x+54=\left(  x-2\right)  \left(  x-3\right)  \left(
x+3\right)  ^{2}.
\]
In this case, $\lambda_{1}=-3,$ $m_{1}=2,\lambda_{2}=3,$ $m_{2}=1$ and $k=4.$
From (\ref{basymp2}), we have%
\begin{align*}
b_{0}(n)  &  \sim-\frac{1}{5}n\left(  -3\right)  ^{n}+\frac{4}{5}\left(
-3\right)  ^{n}-\frac{1}{2}3^{n}\\
b_{1}(n)  &  \sim\frac{1}{10}n\left(  -3\right)  ^{n}-\frac{2}{5}\left(
-3\right)  ^{n}-\frac{1}{12}3^{n}\\
b_{2}(n)  &  \sim\frac{1}{45}n\left(  -3\right)  ^{n}-\frac{4}{45}\left(
-3\right)  ^{n}+\frac{1}{9}3^{n}\\
b_{3}(n)  &  \sim-\frac{1}{90}n\left(  -3\right)  ^{n}+\frac{2}{45}\left(
-3\right)  ^{n}+\frac{1}{36}3^{n}.
\end{align*}
The exact values are%
\begin{align*}
b_{0}(n)  &  =-\frac{1}{5}n\left(  -3\right)  ^{n}+\frac{21}{50}\left(
-3\right)  ^{n}-\frac{1}{2}3^{n}+\frac{27}{25}2^{n}\\
b_{1}(n)  &  =\frac{1}{10}n\left(  -3\right)  ^{n}-\frac{83}{300}\left(
-3\right)  ^{n}-\frac{1}{12}3^{n}+\frac{9}{25}2^{n}\\
b_{2}(n)  &  =\frac{1}{45}n\left(  -3\right)  ^{n}+\frac{2}{225}\left(
-3\right)  ^{n}+\frac{1}{9}3^{n}-\frac{3}{25}2^{n}\\
b_{3}(n)  &  =-\frac{1}{90}n\left(  -3\right)  ^{n}+\frac{11}{900}\left(
-3\right)  ^{n}+\frac{1}{36}3^{n}-\frac{1}{25}2^{n}.
\end{align*}
As we observed before, the main contribution comes from the eigenvalue of
maximum multiplicity, in this case $\lambda_{1}=-3.$
\end{example}

\begin{example}
Finally, let's consider the case of complex eigenvalues of multiplicity
greater than one located on the circle $\left\vert s\right\vert =\rho(A).$ Let%
\begin{align*}
P(x)  &  =x^{5}-9x^{4}+34x^{3}-66x^{2}+65x-25\\
&  =\left(  x-1\right)  \left[  x-\left(  2+\mathrm{i}\right)  \right]
^{2}\left[  x-\left(  2-\mathrm{i}\right)  \right]  ^{2}.
\end{align*}
In this case, $\lambda_{1}=2+\mathrm{i},$ $m_{1}=2,\lambda_{2}=2-\mathrm{i},$
$m_{2}=2$ and $k=5.$ From (\ref{basymp2}), we obtain
\begin{align*}
b_{0}(n)  &  \sim\frac{1}{4}\left(  \sqrt{5}\right)  ^{n}(n-5)\left[
\cos\left(  \theta n\right)  -7\sin\left(  \theta n\right)  \right] \\
b_{1}(n)  &  \sim-\frac{1}{10}\left(  \sqrt{5}\right)  ^{n}(n-5)\left[
2\cos\left(  \theta n\right)  -39\sin\left(  \theta n\right)  \right] \\
b_{2}(n)  &  \sim-\frac{1}{10}\left(  \sqrt{5}\right)  ^{n}(n-5)\left[
2\cos\left(  \theta n\right)  +31\sin\left(  \theta n\right)  \right] \\
b_{3}(n)  &  \sim\frac{1}{10}\left(  \sqrt{5}\right)  ^{n}(n-5)\left[
2\cos\left(  \theta n\right)  +11\sin\left(  \theta n\right)  \right] \\
b_{4}(n)  &  \sim-\frac{1}{20}\left(  \sqrt{5}\right)  ^{n}(n-5)\left[
\cos\left(  \theta n\right)  +3\sin\left(  \theta n\right)  \right]
\end{align*}
with%
\[
\theta=\arctan\left(  \frac{1}{2}\right)  .
\]
The exact values are%
\begin{align*}
b_{0}(n)  &  =\frac{1}{4}\left(  \sqrt{5}\right)  ^{n}\left[  (n-21)\cos
\left(  \theta n\right)  +(-7n+22)\sin\left(  \theta n\right)  \right]
+\frac{25}{4}\\
b_{1}(n)  &  =-\frac{1}{10}\left(  \sqrt{5}\right)  ^{n}\left[  2(n-50)\cos
\left(  \theta n\right)  +(-39n+125)\sin\left(  \theta n\right)  \right]
-10\\
b_{2}(n)  &  =-\frac{1}{10}\left(  \sqrt{5}\right)  ^{n}\left[  (2n+65)\cos
\left(  \theta n\right)  +(31n-100)\sin\left(  \theta n\right)  \right]
+\frac{13}{2}\\
b_{3}(n)  &  =\frac{1}{10}\left(  \sqrt{5}\right)  ^{n}\left[  2(n+10)\cos
\left(  \theta n\right)  +(11n-35)\sin\left(  \theta n\right)  \right]  -2\\
b_{4}(n)  &  =-\frac{1}{20}\left(  \sqrt{5}\right)  ^{n}\left[  (n+5)\cos
\left(  \theta n\right)  +\left(  3n-10\right)  \sin\left(  \theta n\right)
\right]  +\frac{1}{4}.
\end{align*}

\end{example}

\end{document}